\documentclass[conference]{IEEEtran}
\IEEEoverridecommandlockouts
\usepackage{cite}
\usepackage{amsmath,amssymb,amsfonts}
\usepackage{algorithmic}
\usepackage{graphicx}
\usepackage{textcomp}
\usepackage{xcolor}
\usepackage{multirow,bbm,url,mathtools}
\usepackage{graphicx}
\usepackage{caption}
\usepackage{subcaption}
\usepackage{cite}   
\usepackage{verbatim}
\usepackage{wrapfig, bm}
\usepackage{tikz}
\usepackage[ruled,vlined,linesnumbered,algo2e]{algorithm2e}
\usepackage{cuted}
\usepackage{hyperref}
\usepackage{cleveref}
\hypersetup{colorlinks=true,linkcolor=blue,urlcolor=blue,citecolor=blue}
\usepackage{todonotes}
\allowdisplaybreaks
\usetikzlibrary{calc,arrows,fit,decorations.pathreplacing}

\newcommand{\reffig}[1]{Figure~\ref{#1}}
\newcommand{\reftab}[1]{Table~\ref{#1}}

\newcommand{\calN}{\mathcal{N}}
\newcommand{\calG}{\mathcal{G}}
\newcommand{\calL}{\mathcal{L}}
\newcommand{\calD}{\mathcal{D}}
\newcommand{\calT}{\mathcal{T}}

\crefname{equation}{}{}

\def\BibTeX{{\rm B\kern-.05em{\sc i\kern-.025em b}\kern-.08em
    T\kern-.1667em\lower.7ex\hbox{E}\kern-.125emX}}
\begin{document}

\title{
    Scalable Multi-Period AC Optimal Power Flow Utilizing GPUs with High Memory Capacities\thanks{
    Argonne National Laboratory's work was supported by the U.S. Department of Energy, Assistant Secretary for Environmental Management, Office of Science, under contract DE-AC02-06CH11357.
    We gratefully acknowledge the computing resources provided and operated by the Joint Laboratory for System Evaluation (JLSE) at Argonne National Laboratory.
    % {\color{blue} @Michel, which one is the one for GH200?}
  }
} 
% \footnote{\footnotesize The submitted manuscript has been created by UChicago Argonne, LLC, Operator of Argonne National Laboratory (``Argonne''). Argonne, a U.S. Department of Energy Office of Science laboratory, is operated under Contract No. DE-AC02-06CH11357.}

\author
{\IEEEauthorblockN{Sungho Shin, Vishwas Rao, Michel Schanen, D. Adrian Maldonado, and Mihai Anitescu}
  \IEEEauthorblockA{
    \textit{Mathematics and Computer Science Division} \\
    \textit{Argonne National Laboratory}\\
    Lemont, Illinois \\
    \{sshin, vhebbur, mschanen, maldonadod, anitescu\}@anl.gov
  } 
}

\maketitle

\begin{abstract}
  This paper demonstrates the scalability of open-source GPU-accelerated nonlinear programming (NLP) frameworks---ExaModels.jl and MadNLP.jl---for solving multi-period alternating current (AC) optimal power flow (OPF) problems on GPUs with high memory capacities (e.g., NVIDIA GH200 with 480 GB of unified memory).
  There has been a growing interest in solving multi-period AC OPF problems, as the increasingly fluctuating electricity market requires operation planning over multiple periods.
  These problems, formerly deemed intractable, are now becoming technologically feasible to solve thanks to the advent of high-memory GPU hardware and accelerated NLP tools.
  This study evaluates the capability of these tools to tackle previously unsolvable multi-period AC OPF instances. 
  Our numerical experiments, run on an NVIDIA GH200, demonstrate that we can solve a multi-period OPF instance with more than 10 million variables up to $10^{-4}$ precision in less than 10 minutes.
  These results demonstrate the efficacy of the GPU-accelerated NLP frameworks for the solution of extreme-scale multi-period OPF.
  We provide ExaModelsPower.jl, an open-source modeling tool for multi-period AC OPF models for GPUs.
\end{abstract}

\begin{IEEEkeywords}
  power system modeling, optimal power flow, nonlinear programming, GPU computing
\end{IEEEkeywords}

\section{Introduction}

The goal of the AC OPF is to compute the least-cost generation dispatch subject to network physical constraints and security limits. These network constraints are represented in the form of nonlinear power flow equations, resulting in an optimization model that is both nonlinear and nonconvex \cite{frank2016introduction,molzahn2019survey}. To increase the flexibility of the AC OPF, it is possible to endow it with a temporal dimension that will allow modeling processes such as demand and battery response \cite{Capitanescu2016}. This temporal version of the AC OPF, named multi-period AC OPF, consists of a temporally coupled collection of single-period AC OPF problems defined over a time horizon spanning multiple time periods \cite{alizadeh2022envisioning1}. In each period, generator setpoints must be set so that demand in that period is satisfied. The AC OPF models in consecutive time periods are coupled via generator ramping limits.
The multi-period model enables forecasting and planning for the effect of present generator dispatch decisions on future system states and decisions by determining an optimal dispatch policy that is optimized over a desired time horizon. The multi-period formulation offers many advantages over its single-period counterpart, including the ability to make dispatch decisions that account for not only current system conditions (e.g., network loads and topology) but also future forecasted system conditions. For example, our previous work \cite{rao2022frequency} demonstrates the use of multi-period AC OPF for secondary frequency control during contingency events; we also demonstrated that the multi-period formulation can be used to quantify resiliency of the grids \cite{rao2019multiperiod}; it has been used to provision for flexible load modeling in smart grids \cite{usman2021stochastic}; authors of \cite{Jabr2015} demonstrate the use of multi-period formulaton for construction a storage portfolio in presence of renewables in the grid. 

However, while the multiperiod AC OPF provides numerous benefits compared to its traditional counterpart, the increased computational effort of solving this optimization problem has historically limited its application.
Recent advancements in the utilization of GPU accelerators \cite{shinAcceleratingOptimalPower2023, regevHyKKTHybridDirectiterative2023, pacaudAcceleratingCondensedInteriorPoint2023} have now made solving this problem more accessible.
There has been growing interest in utilizing GPU computing resources to solve large-scale NLPs more efficiently \cite{schubigerGPUAccelerationADMM2020,caoAugmentedLagrangianInteriorpoint2016,kimAcceleratedComputationTracking2023,luCuPDLPCStrengthenedImplementation2024,adabagMPCGPURealTimeNonlinear2024,pacaudFeasibleReducedSpace2022,pacaudParallelInteriorPointSolver2023,coleExploitingGPUSIMD2023,pacaudAcceleratingCondensedInteriorPoint2023,shinAcceleratingOptimalPower2023,regevHyKKTHybridDirectiterative2023}.
Traditionally, GPU solvers have relied on factorization-free algorithms, such as first-order methods or second-order methods with iterative linear algebra, e.g., based on preconditioned conjugate gradient \cite{schubigerGPUAccelerationADMM2020,caoAugmentedLagrangianInteriorpoint2016,kimAcceleratedComputationTracking2023,luCuPDLPCStrengthenedImplementation2024,adabagMPCGPURealTimeNonlinear2024,pacaudFeasibleReducedSpace2022}.
While these approaches have achieved some degree of success, ensuring these methods work reliably over a wide range of practical NLPs with a minimal amount of manual parameter tuning has faced challenges.

In the past few years, there has been significant progress in using second-order optimization algorithms with direct linear solvers on GPUs \cite{pacaudFeasibleReducedSpace2022,pacaudParallelInteriorPointSolver2023,coleExploitingGPUSIMD2023,pacaudAcceleratingCondensedInteriorPoint2023,shinAcceleratingOptimalPower2023,regevHyKKTHybridDirectiterative2023}.
These methods are often preferred to factorization-free methods due to their reliable and fast convergence.
With this new technique, solving the static AC OPF for the 78,484-bus system \cite{snodgrassTractableAlgorithmsConstructing} can be performed 10 times faster on GPUs than the state-of-the-art tools running on the contemporary CPUs \cite{shinAcceleratingOptimalPower2023}.
This has been enabled by the advances in the implementation of (i) sparse automatic differentiation (AD), (ii) pivoting-free condensed-space interior point methods (IPMs), and (iii) mature implementation of sparse Cholesky/LDL$^\top$ solvers, such as cuDSS \cite{NVIDIACuDSSPreview}.
Below, we summarize the recent advances in each of these threads.

\paragraph*{Sparse AD}
Evaluating the derivatives of algebraic NLP functions efficiently is crucial in NLP frameworks. This has been particularly important for AC OPF as the model involves algebraic nonlinear power flow equations \cite{coffrinPowerModelsJLOpenSource2018}. While the state-of-the-art AD frameworks, such as JAX, provide GPU-accelerated derivative computation capabilities, they are not optimized for sparse Jacobian and Hessian matrices in large-scale NLPs. Recently, ExaModels.jl,
% \footnote{\url{https://github.com/exanauts/ExaModels.jl}}
a specialized algebraic modeling system (AMS) for classical large-scale NLP, has been developed to facilitate the implementation of GPU-accelerated AD \cite{shinAcceleratingOptimalPower2023}. ExaModels.jl exploits that most large-scale NLPs are expressed with repeated algebraic subexpressions. In addition, the Jacobian and Hessian matrices are highly structured, allowing us to evaluate part of their entries in parallel. These features allow one to reduce the compile-time overhead for applying the chain rule to the subexpressions and to evaluate first and second-order derivatives in parallel. ExaModels.jl provides powerful sparse AD capabilities for GPU and SIMD architectures.

\paragraph*{Pivoting-Free IPMs}
IPM is the state-of-the-art solution algorithm for large-scale constrained NLPs \cite{nocedalNumericalOptimization2006,wachter2006implementation,byrdKnitroIntegratedPackage2006}. Directly running classical IPM on GPUs is considered to be impossible due to the difficulties in parallelizing indefinite LBL$^\top$ factorization of the Karush-Kuhn-Tucker (KKT) system, which requires the numerical pivoting process \cite{swirydowiczGPUresidentSparseDirect2024}.
However, it is possible to parallelize the factorization if the factorization process does not require the use of numerical pivoting.
The pivoting-free solution of KKT systems can be achieved by converting the original indefinite system into a positive definite system, which can be factorized with Cholesky factorization with static pivoting.
Several different approaches have been proposed to allow pivoting-free solutions of KKT systems \cite{regevHyKKTHybridDirectiterative2023, shinAcceleratingOptimalPower2023, pacaudAcceleratingCondensedInteriorPoint2023}.
The hybrid KKT system approach is one of the approaches that enable pivoting-free IPM \cite{regevHyKKTHybridDirectiterative2023}.
In this method, originally proposed by \cite{golubSolvingBlockStructuredIndefinite2003}, the indefinite KKT system is converted to a regularized quasi-definite system, and the KKT system is solved by applying iterative methods to its Schur complement system.
Performing matrix-vector multiplication with the Schur complement requires the use of sparse Cholesky factorization, and parallel factorization methods can be applied here.
Another strategy is based on the lifted KKT system method, where the equality constraints in the original problem are lifted with slack variables and subsequently relaxed to enable condensing the KKT system to a positive definite system \cite{shinAcceleratingOptimalPower2023}.
Lastly, a reduced KKT system strategy is investigated, where the decision variables are partitioned into state and control variables, and the state variables are eliminated to formulate a fully dense, small dimensional KKT system \cite{pacaudAcceleratingCondensedInteriorPoint2023}.
The hybrid KKT system approach was first implemented by a KKT system solver, HyKKT \cite{regevHyKKTHybridDirectiterative2023}, and the lifted KKT system approach was first implemented by MadNLP.jl \cite{shinAcceleratingOptimalPower2023}.
A recent numerical study suggests that the hybrid and lifted KKT system strategy is particularly promising, although the maximum achievable precision is not as high as the CPU counterpart \cite{pacaudApproachesNonlinearProgramming}.
Specifically, the convergence is reliable only up to $\epsilon_{\text{mach}}^{1/4}\approx 10^{-4}$, which is lower than $\epsilon_{\text{mach}}^{1/2}\approx 10^{-8}$ on CPUs (numbers apply for double precision).

\paragraph*{GPU-Accelerated Sparse Linear Solvers}
One of the key requirements of the successful application of the lifted and hybrid KKT system strategies is the availability of an efficient and reliable sparse linear solver that implements either the Cholesky or LDL$^\top$ decompositions.
Recently, there has been significant progress in the implementation efficiency of sparse direct solvers on GPUs.
NVIDIA's CUDA library provides several options, and their new direct solver library, cuDSS, provides impressive performance, as shown in \cite{pacaudGPUacceleratedNonlinearModel2024}.

% \paragraph*{Intro to hardware}
% \begin{itemize}
% \item Recent hardware trend on GPU (unified arch, GPUs within HPC)
% \item Limitations of previous generation GPUs, e.g., multiple GPUs via NVLINK
% \item Their potential significance in optimization
% \end{itemize}
\vspace{-1mm}
\paragraph*{Contributions}
In this paper, we showcase the capabilities of NLP frameworks running on GPUs for solving large-scale multi-period AC OPF problems.
In particular, two recently developed NLP solution tools---ExaModels.jl and MadNLP.jl---provide the capabilities to effectively tackle multi-period AC OPF problems \cite{shinAcceleratingOptimalPower2023}.
The numerical experiments are performed using NVIDIA's GH200 GPU to showcase the modern GPU's capabilities of solving extremely large-scale optimization instances that require using a large amount of memory.
Our results demonstrate that with these software and hardware capabilities, multi-period AC OPF instances with more than 10 million variables can be solved up to $10^{-4}$ accuracy in less than 10 minutes.
We also present an open-source modeling tool for large-scale AC OPF, ExaModelsPower.jl, which provides static and multi-period AC OPF modeling capabilities for GPUs.

\paragraph*{Notations}
We denote the set of real numbers and the set of integers by
$\mathbb{R}$ and $\mathbb{I}$, respectively. We let $[M]:=\{1,2,\cdots,M\}$.
A vector of ones with an
appropriate size is denoted by $\boldsymbol{1}$. An identity matrix
with an appropriate size is denoted by $I$. For matrices $A$ and $B$,
$A\succ(\succeq) B$ indicates that $A-B$ is positive (semi-)definite
while $A>(\geq)B$ denotes a componentwise inequality. We use the
convention $X:=\mathop{\text{diag}}(x)$ for any symbol $x$.

\section{Multi-period AC OPF}
% \begin{itemize}
% \item AC OPF formulation
% \item Multi-period formulations (ramping constraints)
% \end{itemize}
The goal of the classical single-period AC OPF is to find an optimal generation dispatch, along with the network voltage values, such that all network quantities (e.g., real and reactive power generation, bus voltages, and transmission line currents) stay within allowed physical limits, and the resulting network power flows satisfy Kirchhoff’s laws. These physical laws are represented in the form of power flow equations, which are nonlinear, and the resulting model is typically both nonlinear and nonconvex \cite{frank2016introduction,molzahn2019survey}. The multi-period AC OPF model is a temporally coupled collection of single-period AC OPF problems defined over a time horizon spanning multiple time periods. Each period consists of a classical AC OPF model where generator setpoints must be set so that all electrical loads in that period are satisfied. 

The multi-period OPF considered in this work is formulated as follows:
\begin{subequations}
    \label{eq:opf}
    \begin{align}
      \min \;
      & \sum_{g \in \calG} \sum_{t\in\calT} f_g\left(p^G_{gt}\right) \label{opf:obj}\\
      \text{s.t.} \;
      & \sum_{l\in\calL_n^\text{to}} p_{lt} - \sum_{l\in\calL_n^\text{from}} p_{lt} + \sum_{g\in\calG_n} p_{gt}^G = \sum_{j\in\calD_n} P_{jt}^D \notag \\
      & \quad \forall n\in\calN, \; t\in\calT, \label{eq:balancep}\\
      & \sum_{l\in\calL_n^\text{to}} q_{lt} - \sum_{l\in\calL_n^\text{from}} q_{lt} + \sum_{g\in\calG_n} q_{gt}^G = \sum_{j\in\calD_n} Q_{jt}^D \notag \\
      & \quad \forall n\in\calN, t\in\calT, \label{eq:balanceq}\\
      & p_{lt} = v_{mt} v_{nt} \left[ G_l \cos(\theta_{nt} - \theta_{mt}) + B_l \sin(\theta_{nt} - \theta_{mt}) \right] \notag \\
      & \quad \forall l=(m,n)\in\calL, \; t\in\calT, \label{eq:pfp}\\
      & q_{lt} = v_{mt} v_{nt} \left[ G_l \sin(\theta_{nt} - \theta_{mt}) - B_l \cos(\theta_{nt} - \theta_{mt}) \right] \notag \\
      & \quad \forall l=(m,n)\in\calL, \; t\in\calT, \label{eq:pfq}\\
      %& | q_{gt}^G - q_{g,t-1}^G | \leq Q_g^\text{ramp} \quad \forall g\in\calG, \; t\in\calT\backslash\{1\}, \label{eq:rampq}\\
      & \sqrt{p_{lt}^2 + q_{lt}^2} \leq S_l^\text{max} \, \quad \forall l\in\calL, \; t\in\calT,\label{eq:linecap}\\
      & P^\text{min}_g \leq p^G_{gt} \leq P^\text{max}_g\, \quad \forall g\in\calG, \; t\in\calT, \label{opf:pgencon}\\
      & Q^\text{min}_g \leq q^G_{gt} \leq Q^\text{max}_g\, \quad \forall g\in\calG, \; t\in\calT, \label{opf:qgencon} \\
      & V_n^\text{min} \leq v_{nt} \leq V_n^\text{max} \quad \forall n\in\calN, \; t\in\calT, \label{opf:vmagcon} \\
      & \Theta_{mn}^\text{min} \leq \theta_{mt} - \theta_{nt} \leq \Theta_{mn}^\text{max} \quad \forall (m,n)\in\calL, \; t\in\calT. \label{eq:anglerange}\\
      & | p_{gt}^G - p_{g,t-1}^G | \leq P_g^\text{ramp} \quad \forall g\in\calG, \; t\in\calT\backslash\{1\}, \label{eq:rampp}
    \end{align}
    \end{subequations}
    % where the objective function represents the total generation cost of the system.
    Here, $\calN$ represents the set of buses, $\calD_n$ represents the set of loads at bus $n\in\calN$, $\calG_n$ denotes the set of generators at bus $n\in\calN$, $\calL$ is used to denote set of transmission lines, $\calL_n^\text{from}$ represent the set of transmission lines connected from bus $n$, $\calL_n^\text{to}$ represents the set of transmission lines connected to bus $n$, and $\calT$ the set of time periods. Furthermore, $G_l$ denotes the conductance on line $l$, $B_l$ susceptance on line $l$, $P_j^D$ denotes real power required at load $j$, $Q_j^D$ denotes reactive power required at load $j$, $S_l^\text{max}$ denotes apparent power capacity on line $l$, $P_g^\text{min}$ denotes minimum real power capacity at generator $g$, $P_g^\text{max}$ denotes maximum real power capacity at generator $g$, $P_g^\text{ramp}$ denotes real power ramp capacity at generator $g$, $Q_g^\text{min}$ denotes minimum reactive power capacity at generator $g$, $Q_g^\text{max}$ denotes maximum reactive power capacity at generator $g$, $Q_g^\text{ramp}$ denotes reactive power ramp capacity at generator $g$, $V_n^\text{min}$ denotes minimum voltage magnitude at bus $n$, $V_n^\text{max}$ denotes maximum voltage magnitude at bus $n$, $\Theta_n^\text{min}$ denotes minimum voltage angle at bus $n$, $\Theta_n^\text{max}$ denotes maximum voltage angle at bus $n$, $p_{gt}^G$ denotes real power production at generator $g$ at time $t$, $q_{gt}^G$ denotes reactive power production at generator $g$ at time $t$, $p_{lt}$ denotes real power flow on line $l$ at time $t$, $q_{lt}$ denotes reactive power flow on line $l\in\calL$ at time $t\in\calT$, $v_{nt}$ denotes voltage magnitude at bus $n\in\calN$ at time $t\in\calT$, and $\theta_{nt}$ denotes phase angle at bus $n\in\calN$ at time $t\in\calT$.
    
    Constraints~\eqref{eq:balancep} and \eqref{eq:balanceq} balance the power generation and load for each bus $n\in\calN$. Constraints~\eqref{eq:pfp} and \eqref{eq:pfq} represent the power flow equations from Kirchoff's law. The other bound constraints~\eqref{eq:linecap} -- \eqref{eq:anglerange} formulate the physical capacities of the system. Note that the set of equations in \eqref{eq:opf} reduces to the single-period AC OPF by setting $\mathcal{T} =\{1\}$.  The ramping of power production is restricted by generator capacity in constraint~\eqref{eq:rampp}. The coupling between different periods arises from these constraints. More details of the static OPF formulations can be found in \cite{coffrinPowerModelsJLOpenSource2018}.

\section{NLP Frameworks for GPUs}
We now describe the key ideas in the recent advances in NLP software frameworks for GPUs.
In particular, we discuss ExaModels.jl, an algebraic modeling and AD tool, and MadNLP.jl, an NLP solver.
% While sparse linear algebra, in particular, Cholesky or LDL$^\top$ factorization, plays important role, we do not discuss it here, as we are relying on NVIDIA's proprietary library.

\subsection{Sparse AD}
AD is the most preferred way of evaluating derivatives of computer programs due to their superior accuracy, efficiency, and ease of use compared to their alternatives.
AD tools directly convert the function evaluation code into a derivative evaluation code.
There are two different ways to apply recursion:  forward mode and reverse mode.
Reverse-mode AD (adjoint method) has proven to be particularly effective in dealing with function expressions in large-scale optimization problems.
There exist several different ways of implementing AD, including operator overloading and source code transformation, both of which can achieve high performance.

When it comes to large-scale NLP, such as AC OPFs, AD capabilities are often provided by AMSs, like JuMP \cite{dunningJuMPModelingLanguage2017}, CasADi \cite{anderssonCasADiSoftwareFramework2019}, and AMPL \cite{fourerModelingLanguageMathematical1990}.
These tools are more specialized for the evaluation of sparse Jacobian and Hessian matrices encountered within large-scale optimization.
ExaModels.jl can also be classified as an AMS for mathematical programs.
ExaModels.jl is built on Julia Language's multiple dispatch feature, which allows the function to be dispatched based on the run-time types of the argument \cite{bezansonJuliaFreshApproach2017}.
This feature enables straightforward implementation of operator overloading-based automatic differentiation.
ExaModels.jl is different from the other AD tools in that they are designed to {\it exploit the repeated patterns} within the NLP models.
A similar idea has been implemented in one of the recently developed AMS, Gravity \cite{hijaziGravityMathematicalModeling2018a}, and this feature is called constraint template.

ExaModels.jl exploits the repeated expression structures by storing NLP problem data in a structured format by leveraging what's called {\it SIMD abstraction of NLPs}. This abstraction can be mathematically expressed as follows:
\begin{align}\label{eqn:prob}
  \min_{x^\flat\leq x \leq x^\sharp}
  & \sum_{l\in[L]}\sum_{i\in [I_l]} f^{(l)}(x; p^{(l)}_i)\\\nonumber
  \text{s.t.}\; &\sum_{n\in [N_m]}\sum_{k\in [K_n]}h^{(n)}(x; s^{(n)}_{k}) =0,\quad \forall m\in[M]
\end{align}
where $f^{(\ell)}(\cdot,\cdot)$ and $h^{(n)}(\cdot,\cdot)$ are twice differentiable functions with respect to the first argument, whereas $\{\{p^{(k)}_i\}_{i\in [N_k]}\}_{k\in[K]}$ and $\{\{\{s^{(n)}_{k}\}_{k\in[K_n]}\}_{n\in[N_m]}\}_{m\in[M]}$ are problem data.
One can observe that the problem in \cref{eqn:prob} is expressed with repeated expressions, which can be evaluated and differentiated in parallel.
Accordingly, the evaluation and differentiation of the model equations in \cref{eqn:prob} are amenable to SIMD parallelism on GPU hardware.
The user interface in ExaModels.jl mandates that users define their model equations using a {\tt Generator} data type.
{\tt Generator} is a composite data type composed of a command (a function in Julia) and data (an array located on either the host or a computing device) on which the command operates.
This abstraction enables ExaModels.jl to store NLP data in a structured format, which in turn facilitates performing AD on GPUs.

The standard reverse-mode AD can be applied to the functions in \cref{eqn:prob}.
Here, the key novelty of ExaModels.jl's implementation is that reverse-mode AD is applied to each pattern.
This allows obtaining compiled AD kernels that can be executed over multiple data.
These compiled AD kernels are highly efficient, as they are specifically compiled for each computational pattern.
Also, these AD kernels can be run over data on various GPU array formats.
In turn, ExaModels.jl can provide the capability to evaluate efficiently first- and second-order derivatives on diverse architectures of GPU accelerators.

Typically, sparse AD relies on matrix coloring and Jacobian/Hessian-vector product evaluations.
However, ExaModels.jl directly evaluates sparse Jacobian and Hessian matrices in coordinate list (COO) format with pre-analyzed sparsity patterns.
We perform sparsity analysis for each of the patterns and construct the sparse Jacobian/Hessian matrices by expanding the analyzed sparsity pattern over the data array.
Optionally, these sparse matrices can be compressed further to compressed sparse column (CSC) formats if desired.
For a more detailed description of the AD implementation, the readers are referred to our previous publication \cite{shinAcceleratingOptimalPower2023}.

\subsection{Condensed-Space IPMs: Lifted KKT Approach}

Among different approaches of pivoting-free IPM, here we explain the lifted KKT system approach. For hybrid KKT system or reduction approaches, the readers are referred to \cite{regevHyKKTHybridDirectiterative2023} and \cite{pacaudAcceleratingCondensedInteriorPoint2023}, respectively.

Consider the following standard form NLP:
\begin{equation}\label{eqn:cpt}
    \min_{x^\flat \leq x \leq x^\sharp}\;  f(x) \quad \text{s.t.}\;
     g(x) =0 .
\end{equation}
In the classical implementation of IPM,
the following indefinite system, the so-called KKT system, is
solved using sparse linear solvers to compute the step direction, typically based on
LBL$^\top$ factorization:
\begin{align}\label{eqn:kkt-indefinite}
  &\begin{bmatrix}
    W^{(\ell)}  + \Sigma^{(\ell)} + \delta^{(\ell)}_w I& A^{(\ell)\top}\\
    A^{(\ell)} & -\delta_c^{(\ell)} I\\
  \end{bmatrix}
  \begin{bmatrix}
    \Delta x\\
    \Delta y\\
  \end{bmatrix}=
  \begin{bmatrix}
    r_x^{(\ell)}\\
    r_y^{(\ell)}\\
  \end{bmatrix},
\end{align}
where
\begin{align*}
  W^{(\ell)}
  &:=\nabla^{2}_{xx}\mathcal{L}(x^{(\ell)},y^{(\ell)},z^{\flat(\ell)},z^{\sharp(\ell)}),\quad
  A^{(\ell)}
  := \nabla_xg(x^{(\ell)})\\
  \Sigma^{(\ell)}&:= (X^{(\ell)})^{-1}Z^{(\ell)},\quad
  r_x^{(\ell)}
  :=\nabla_x f(x^{(\ell)}) - \mu (X^{(\ell)})^{-1} \boldsymbol{1},\\
  r_y^{(\ell)}
  &:=g(x^{(\ell)}),
\end{align*}
and $\delta^{(\ell)}_w, \delta^{(\ell)}_c>0$ are determined through the inertia correction procedure \cite{wachterImplementationInteriorpointFilter2006}.
Line search is employed to determine the step size associated with the step direction determined by \cref{eqn:kkt-indefinite}. 

The lifted KKT method has two key differences from the standard IPM procedure: (i) inequality relaxation and (ii) condensation of the KKT system.
We first describe the inequality relaxation procedure.
In this step, we relax the original equality constraints in \cref{eqn:cpt} by introducing slack variables $s\in\mathbb{R}^{m}$ and allow them to vary slightly around zero:
\begin{align}\label{eqn:relax}
  g(x)- s = 0,\quad s^{\flat}\leq s\leq  s^\sharp,
\end{align}
where $s^\flat,s^\sharp\in\mathbb{R}^{m}$ are lower and upper bounds chosen by the solver. In our implementation in MadNLP.jl, we choose them to be $\pm \epsilon_{\text{tol}} \boldsymbol{1}$, where $\epsilon_{\text{tol}}>0$ is a user-specified relative tolerance
of the IPM.
This setting ensures the constraint violation within the desired relative error bound.
% The inequality relaxation procedure yields the following problem:
% \begin{equation}\label{eqn:relaxed}
%   \min_{\big[\substack{x^\flat\\s^\flat}\big]\leq \big[\substack{x\\s}\big] \leq \big[\substack{x^\sharp\\s^\sharp}\big]}\;
%   f(x)\quad \text{s.t.}\;  g(x) - s = 0.
% \end{equation}
With the additional structures introduced by \cref{eqn:relax}, one can write the KKT system \cref{eqn:kkt-indefinite} as

\begin{strip}
  \rule{\textwidth}{.5pt}
  \begin{align}
    \underbrace{
    \begin{bmatrix} 
      W^{(\ell)}  + \delta^{(\ell)}_w I & & A^{(\ell)\top}& -I & I &  \\
                                        & \delta^{(\ell)}_w I & -I&&&-I & I\\
      A^{(\ell)}& -I & -\delta^{(\ell)}_c I\\
      Z_x^{(\ell)\flat}&&&X^{(\ell)}-X^\flat\\
      -Z_x^{(\ell)\sharp}&&&&X^\sharp-X^{(\ell)}\\
                                        &Z_s^{(\ell)\flat}&&&&S^{(\ell)}-S^\flat\\
                                        &-Z_s^{(\ell)\sharp}&&&&&S^\sharp-S^{(\ell)}\\
    \end{bmatrix}
    }_{M_\text{full}}
    \begin{bmatrix}
      \Delta x \\
      \Delta s \\
      \Delta y \\
      \Delta z_x^\flat \\
      \Delta z_x^\sharp \\
      \Delta z_s^\flat \\
      \Delta z_s^\sharp \\
    \end{bmatrix} =
    \begin{bmatrix}
      p^{(\ell)}_{x }\\
      p^{(\ell)}_{s }\\
      p^{(\ell)}_{y }\\
      p^{(\ell)}_{z_x^\flat }\\
      p^{(\ell)}_{z_x^\sharp }\\
      p^{(\ell)}_{z_s^\flat }\\
      p^{(\ell)}_{z_s^\sharp }\\
    \end{bmatrix},
    % \tag{Full KKT System}
    \label{eqn:very-long-eqn} 
  \end{align}
  \rule{\textwidth}{.5pt} 
\end{strip}
\noindent where $p^{(\ell)}_x,\cdots p^{(\ell)}_{z_s^\sharp}$ are defined as follows:
\begin{align*}
  p_x&=\nabla_{x} f(x) - \nabla_{x}g(x)^\top y  - z_x^\flat  + z_x^\sharp \\
  p_s&=- z_s^\flat  + z_s^\sharp + y ,\quad
  p_y=g(x) - s  \\
  p_{z^\flat_x}&=Z^\flat_x (x-x^\flat) - \mu\boldsymbol{1} ,
\quad  p_{z^\sharp_x}=Z^\sharp_s (s-s^\sharp) - \mu\boldsymbol{1},\\
  p_{z^\flat_s}&=Z^\sharp_x (x^\sharp-x) - \mu\boldsymbol{1}
  \quad p_{z^\sharp_S}=Z^\sharp_s (s^\sharp-s) - \mu\boldsymbol{1}.
\end{align*}
Here, the iteration index $\ell$ is dropped for simplicity.
% $p^{(\ell)}_x,\cdots p^{(\ell)}_{z_s^\sharp}$ are defined by the left-hand sides of the equations in
% \cref{eqn:first}.

The primary benefit of the inequality relaxation is the introduction of extra structure to the enhanced KKT system.
Specifically, the lower-right 6x6 block in \cref{eqn:very-long-eqn} can be
eliminated due to the invertibility of the lower-right, which derives from $\delta_w,\delta_c, x^{(\ell)}-x^\flat, x^\sharp-x^{(\ell)}, z^{\flat(\ell)}, z^{\sharp(\ell)}> 0$.
As a result of this elimination, often referred to as {\it condensation}, we can obtain the following {\it condensed KKT system}:
\begin{align*}
 (\underbrace{W + \delta_wI + \Sigma_x + A^{\top} D A}_{M_\text{cond}} ) \Delta x = q_x + A^\top (C q_s +  Dq_y ),
\end{align*}
where
\begin{align*}
  \Sigma_x&:= Z^\flat_x (X-X^\flat)^{-1}+ Z^\sharp_x (X^\sharp-X)^{-1}\\
  \Sigma_s&:= Z^\flat_s (S-S^\flat)^{-1}+ Z^\sharp_s (S^\sharp-S)^{-1}\\
  q_x&:=p_x + (X-X^\flat)^{-1} p_{z^\flat_x}-  (X^\sharp-X)^{-1} p_{z^\sharp_x}\\
  q_s&:=p_s + (S-S^\flat)^{-1} p_{z^\flat_s}-  (S^\sharp-S)^{-1} p_{z^\sharp_s}\\
  q_y&:=p_y\\
  C &:= \left(\delta_c \Sigma_s + (1+\delta^{}_c\delta^{}_w) I\right)^{-1}, \quad
  D := \left(\Sigma_s + \delta^{}_w I\right)C.
\end{align*}
Once the condensed system is solved, the dual and slack step directions can be recovered straightforwardly with:
\begin{align*}
  \Delta s &:= C \left(\delta_c q_s - (q_y + A\Delta x)\right)\\
  \Delta y &:= (\Sigma_s + \delta_w I) \Delta s -q_s.\label{eqn:recover-2}\\
  \Delta z^\flat_x &= \left(X-X^\flat\right)^{-1} \left(-Z^\flat_x \Delta x  + p_{z^\flat_x}\right)\\
  \Delta z^\sharp_x &= \left(X^\sharp-X\right)^{-1} \left(Z^\sharp_x \Delta x  + p_{z^\sharp_x}\right)\\
  \Delta z^\flat_s &= \left(S-S^\flat\right)^{-1} \left(-Z^\flat_s \Delta s  + p_{z^\flat_s}\right)\\
  \Delta z^\sharp_s &= \left(S^\sharp-S\right)^{-1} \left(Z^\sharp_s \Delta s  + p_{z^\sharp_s}\right).
\end{align*}
It is important to note that the only required sparse matrix factorization is for $M_\text{cond}$.
Upon the application of standard inertia correction procedures, it is guaranteed that $M_\text{cond} \succ 0$ (see the discussion in \cite{shinAcceleratingOptimalPower2023}).
Thus, Cholesky or LDL$^\top$ factorization with static pivoting can be applied to solve the condensed system.
Consequently, this helps us to overcome the difficulties in numerical pivoting and makes the implementation of IPM on GPUs possible.
For a more detailed description of the condensed-space IPM implementation, the readers are referred to our previous publication \cite{shinAcceleratingOptimalPower2023}.

\section{Numerical Results}
We now present the numerical results demonstrating the NLP solution capabilities on GPUs. The hardware and software details for the configurations can be found in \Cref{tab:systems}. The multi-period AC OPF model is formulated with our new open-source modeling tool ExaModelsPower.jl.\footnote{\url{https://github.com/exanauts/ExaModelsPower.jl}} The script to reproduce the numerical results is available online.\footnote{\url{https://github.com/exanauts/multiperiod-osmses-2024}}
% More detailed version and dependency information can be found in the {\tt Manifest.toml} file in the script repository.

\begin{table*}[t]
  \centering
  \caption{Hardware and software configurations}
  \begin{tabular}{|c|c|c|c|}
    \hline
    Label & AD Tool & NLP Solver & Linear Sovler\\
    \hline
    AMD EPYC 7532& ExaModels.jl v0.7.1 & Ipopt v3.14.14 & Ma27 (libHSL-2023.11.7)\\
    NVIDIA A100 40GB & ExaModels.jl v0.7.1 & MadNLP v0.8.3 & cuDSS v0.2.1\\
    NVIDIA GH200 480GB & ExaModels.jl v0.7.1 & MadNLP v0.8.3 & cuDSS v0.2.1\\
    \hline
  \end{tabular}
  \label{tab:systems}
\end{table*}

The systems used are workstations equipped with the CPU and GPUs listed in \reftab{tab:systems}. The CPU is an EPYC 7532 with 32 cores, although our linear solver of choice (HSL MA27) does not use multi-threading. The GPUs of choice are an NVIDIA A100 with 32 GB of RAM and an NVIDIA GH200 with 480GB of unified RAM with 96 GB of HBM3 memory. The GH200 is particularly interesting for multiperiod models as it provides the highest amount of RAM through the unified RAM setup. On this machine, we allocated all memory as unified memory arrays. We expect such a setup of multiple GPUs connected by unified RAM through a fast bus or network to become commodity hardware in the near future.

% All benchmarks were run with the same software configuration except for the linear solver. On the CPU, we used MA27 of the linear solver collection HSL, and on the GPU, we used the latest available version of cuDSS, a new experimental sparse linear solver by NVIDIA.

\begin{table}
\centering
\caption{List of cases benchmarked}
\begin{tabular}{|l|c|c|c|c|}
\hline 
{\bf Case}  & {\bf resolution} & {\bf nperiods}& {\bf nvars} & {\bf ncons}\\
% \multirow{2}{*}{\bf case} & \multirow{2}{*}{\bf nvars} & \multirow{2}{*}{\bf ncons} & \multicolumn{2}{c|}{\textbf{GPU}} & \multicolumn{2}{c|}{\textbf{GPU}} \\
% \cline{4-7} &  &  & iter & solution time (sec) & iter & solution time (sec) \\
\hline 
case30 & 30 min & 30 &   7k &  11k \\ 
case118& 60 min &  168 &  183k & 268k \\
case118 & 5 min &   10080 & 11m &  16m \\
case1354 & 30 min &  30 & 336k & 507k \\
case9241 & 30 min & 30 &  3m &   4m \\
case9241 & 30 min & 180 & 15m &  24m \\
\hline 
\end{tabular}
\label{tab:cases}
\end{table}

\reftab{tab:cases} lists all the cases we used in our benchmarks. The base cases are IEEE-30, IEEE-118, IEEE-1354, and PEGASE-9241 with 30, 118, 1354, and 9241 buses, respectively. The loads for the periods were generated as a random perturbation of a scaled periodic function (such as a sinusoid) to reflect the fluctuations in demand across time. The exact datasets used for the experiments will be made available in the git repository. %{\color{blue}@Vishwas could you explain this?} %\todo{How were the loads generated?}.

\begin{figure}
    \centering
    \includegraphics[width=0.47\textwidth]{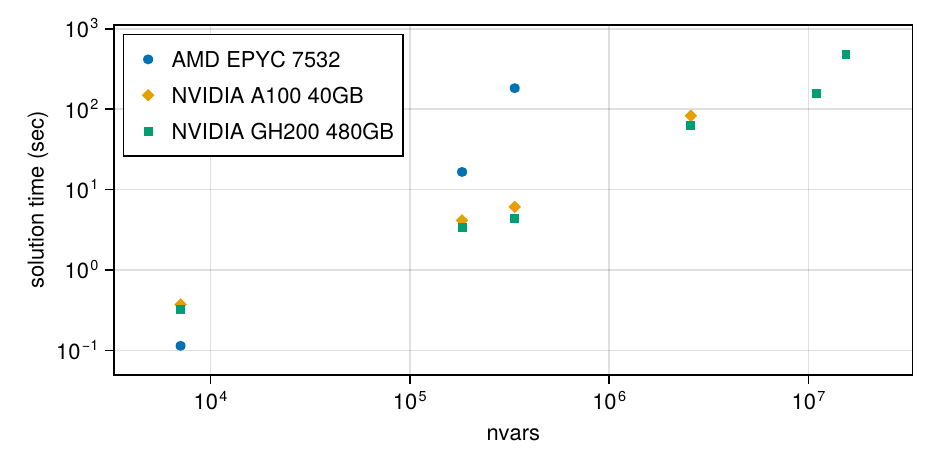}
    \caption{Benchmark results for the cases listed in \reftab{tab:cases}.}
    \label{fig:results}
\end{figure}

In \reffig{fig:results}, we show the benchmark results of all cases listed in \reftab{tab:cases} on the systems in \reftab{tab:systems}.  The two largest cases IEEE-118 with 100,080 periods (10m vars and PEGASE9241 with 180 periods (15m vars) were only able to complete on the GH200, because the A100 ran out of memory. The sequential linear solver MA27 does not provide any parallelization through multithreading and is thus unable to scale with the number of variables. On the other hand, the GPUs can leverage their many-core architecture to parallelize both the function evaluations and linear solver. cuDSS is a closed-source solver and we have no insight into the details of the algorithm. However, a profiling result using NVIDIA NSight Systems shows transfers from the GPU to the CPU on the A100/EPYC system and CPU usage on both the A100 and GH200. This might explain the slightly better performance of the GH200 compared to the A100, as the GH200 relies on uniform memory and does not require slow data transfers between GPU and system RAM. The GH200 is the only choice for the two largest cases, IEEE-118 with 100,080 periods (10m vars) and PEGASE9241 with 180 periods (15m vars), as these problems do not fit into the memory of a single A100 GPU.

\section{Conclusions and Future Outlook}
% write a conclusion based on 
%   This paper showcases the capabilities of the graphics processing unit (GPU)-accelerated nonlinear programming (NLP) frameworks---ExaModels.jl and MadNLP.jl---for solving extreme-scale instances of multi-period alternating current (AC) optimal power flow (OPF) on GPUs with high memorys.
%   There has been growing interest in solving multi-period AC OPF problems, as the increasingly fluctuating electricity market requires operation planning over multiple periods.
%   These problems, formerly deemed intractable, are now becoming technologically feasible to solve thanks to the advent of high-memory GPU hardware and accelerated NLP tools. This study aims to investigate these tools' efficacy in addressing multi-period AC OPF instances that were previously unsolvable. Our numerical experiments, run on an NVIDIA GH200 GPU with 480GB of memory, demonstrate that we can solve a multi-period OPF instance with more than 10 million variables up to $10^{-4}$ precision in less than 3 minutes.
%   These results demonstrate the efficacy of the GPU-accelerated NLP frameworks for the solution of extreme-scale multi-period OPF.
%   We provide ExaModelsPower.jl, a reference implementation of static and multi-period AC OPF models for GPUs.

We have demonstrated the potential of GPU-accelerated NLP frameworks, particularly ExaModels.jl, and MadNLP.jl, in solving large-scale multi-period AC OPF problems that were once considered infeasible due to computational constraints. The previously reported NLP framework running on GPUs has been further strengthened by the use of GPUs with high-capacity memory. With the enhanced hardware capabilities, we have demonstrated that multi-period AC OPF instances with more than 10 million variables can be solved in less than 10 minutes. We expect CPU/GPU systems with unified memory, like the NVIDIA GH200, will become mainstream hardware in the near future with even higher RAM capacities. These advances can provide powerful computational tools for electricity market operations.
We have presented ExaModelsPower.jl as a reference implementation for the multi-period AC OPF model for GPUs. In the future, we plan to further extend the modeling capabilities to cover other large-scale formulations, which may also benefit from the utilization of computational power of GPUs, such as security-constrained AC OPF, frequency control problems, and problems with energy storage.
\bibliographystyle{IEEEtran}
\bibliography{refs} 

 \noindent\fbox{\parbox{\linewidth}{\footnotesize
     The submitted manuscript has been created by UChicago Argonne, LLC, Operator of Argonne National Laboratory (``Argonne''). Argonne, a U.S. Department of Energy Office of Science laboratory, is operated under Contract No. DE-AC02-06CH11357. The U.S. Government retains for itself, and others acting on its behalf, a paid-up nonexclusive, irrevocable worldwide license in said article to reproduce, prepare derivative works, distribute copies to the public, and perform publicly and display publicly, by or on behalf of the Government. The Department of Energy will provide public access to these results of federally sponsored research in accordance with the DOE Public Access Plan (http://energy.gov/downloads/doe-public-access-plan).}
 }
\end{document}